\documentclass[11pt]{article}
\usepackage[left=1in,right=1in,top=1in,bottom=1in]{geometry}
\usepackage{times}
\usepackage{expl3}
\usepackage{cite}
\usepackage[table]{xcolor}
\usepackage{multirow}
\usepackage{stackengine} 
\usepackage{hhline}
\usepackage{lipsum}
\usepackage{titlesec}
\usepackage{wrapfig}
\usepackage{epsfig}
\usepackage{graphicx}
\usepackage{amsmath}
\usepackage[title]{appendix}
\usepackage{amssymb}
\usepackage{epstopdf}
\usepackage{boldline}
\usepackage{calligra}
\usepackage{bm}
\usepackage{url}
\usepackage{blindtext}

\newcommand{\define}{\stackrel{\mbox{\tiny def}}{=}}

\newtheorem{definition}{Definition}
\newtheorem{theorem}{Theorem}
\newtheorem{corollary}{Corollary}
\newtheorem{lemma}{Lemma}

\usepackage{mathtools}
\usepackage{epstopdf}
\usepackage{balance}
\usepackage{thmtools}
\usepackage{thm-restate}
\usepackage{hyperref}
\usepackage{cleveref}

\usepackage[ruled,vlined]{algorithm2e}
\include{pythonlisting}
\newcommand{\ostar}{\mathbin{\mathpalette\make@circled\star}}

\newcommand*{\Scale}[2][4]{\scalebox{#1}{$#2$}}%
\makeatletter
\newcommand{\removelatexerror}{\let\@latex@error\@gobble}
\makeatother
\makeatletter
\newcommand*{\rom}[1]{\expandafter\@slowromancap\romannumeral #1@}
\makeatother

\ExplSyntaxOn
\newcommand\latinabbrev[1]{
  \peek_meaning:NTF . {% Same as \@ifnextchar
    #1\@}%
  { \peek_catcode:NTF a {% Check whether next char has same catcode as \'a, i.e., is a letter
      #1.\@ }%
    {#1.\@}}}
\ExplSyntaxOff

\titleclass{\subsubsubsection}{straight}[\subsubsection]

\begin{document}
\vspace{1cm}
\title{Tail Bounds for Tensor-valued Random Process}\vspace{1.8cm}
\author{Shih~Yu~Chang 
% <-this % stops a space
\thanks{Shih Yu Chang is with the Department of Applied Data Science,
San Jose State University, San Jose, CA, U. S. A. (e-mail: {\tt
shihyu.chang@sjsu.edu}).
           }}

\maketitle

\begin{abstract}
To consider a high-dimensional random process, we propose a notion about stochastic tensor-valued random process (TRP).  In this work, we first attempt to apply a generic chaining method to derive tail bounds for all $p$-th moments of the supremum of TRPs. We first establish tail bounds for suprema of processes with an exponential tail, and further derive tail bounds for suprema of processes with arbitrary number of exponential tails. We apply these bounds to high-dimensional compressed sensing and empirical process characterizations.  
\end{abstract}


\section{Introduction}\label{sec:Introduction} 

% GC

Generic chaining is a powerful probabilistics tool developed by Talagrand to determine the expected value of the supremum of a real-valued stochastic process. This approach stems from the traditional chaining argument invented by Kolmogorov and the majorizing measures method invented by Dudley, Fernique and Talagrand~\cite{talagrand2005generic}. Generic chaining applies $\gamma$-functionals, a quantitative measure of the metric complexity of the index set of the process, to estimate the expectation of the supremum of a random process. Such estimations are known to be sharp at several special cases, e.g., the celebrated majorizing measure theorem provides the sharpness for Gaussian random processes~\cite{talagrand1987regularity}. The applications of generic chaining can be discovered at signal processing, statistics, and harmonic analysis~\cite{talagrand2005generic}. Besides knowing the upper bound for the expected supremum of a random process, we wish to know how probable it is that the supremum of the random process exceeds the expectation of the upper bound. To answer this question, a generic chaining bound has to provide other tail bounds for the deviation of the supremum with respect to the expected supremum of a random process. Such concentration studies of random variables have generated many research results about tail bounds estimations and produced various applications in various fields, see~\cite{ledoux2001concentration}.

Tensors, high-dimensional data, have been applied to different domains in science and engineering, e.g., theoretical physics~\cite{hess2015tensors}, signal processing~\cite{chang2022TWF}, machine learning~\cite{ji2019survey}, etc. Concentration of random matrices (order 2 tensors) has begun its study since 2010~\cite{tropp2012user}, and the generalization of concentration theory to random tensors with arbitrary order can be found at following works: non-independent random tensors by expander mathod~\cite{chang2021tensorExp}, majorization approach~\cite{chang2022general}, T-product tensors~\cite{chang2022T_I,chang2022T_II}, and others~\cite{chang2021convenient,chang2022generalized}. In this work, we first attempt to apply a generic chaining method to derive tail bounds for all $p$-th moments of the supremum of a stochastic tensor-valued random process (TRP). Our work follows the same generic chaining approach given by~\cite{dirksen2015tail}, however, we extend this work to random tensors settings, instead of random variables, and consider the mixed TRP with an arbitrary number of exponential tails with different exponents.

The rest of this paper is organized as follows. In Section~\ref{sec:Tensor-valued Random Process}, we will define tensor-valued random process (TRP) and review notions about tensors briefly. We then consider two types of tail bounds, exponential tail and mixed tail, based on generic chaining in Section~\ref{sec:Supremum and Tail Estimates for TRP} . The application of exponential tail bound to compressed sensing is discussed in Section~\ref{sec:Application: Compressed Sensing}. The application of mixed tail bound to the empirical process is presented in Section~\ref{sec:Application: Empirical Process}.

\section{Tensor-valued Random Process}\label{sec:Tensor-valued Random Process}

We will define tensors considered in this work and tensor-valued random process. More details about basic tensor concepts can be found in~\cite{chang2021convenient,chang2022TWF}. Without loss of generality, one can partition the dimensions of a tensor into two groups, say $M$ and $N$ dimensions, separately. Thus, for two order-($M$+$N$) tensors: $\mathcal{A} \define (a_{i_1, \cdots, i_M, j_1, \cdots,j_N}) \in \mathbb{C}^{I_1 \times \cdots \times I_M\times
J_1 \times \cdots \times J_N}$ and $\mathcal{B} \define (b_{i_1, \cdots, i_M, j_1, \cdots,j_N}) \in \mathbb{C}^{I_1 \times \cdots \times I_M\times
J_1 \times \cdots \times J_N}$, according to~\cite{liang2019further}, the \emph{tensor addition} $\mathcal{A} + \mathcal{B}\in \mathbb{C}^{I_1 \times \cdots \times I_M\times
J_1 \times \cdots \times J_N}$ is given by 
\begin{eqnarray}\label{eq: tensor addition definition}
(\mathcal{A} + \mathcal{B} )_{i_1, \cdots, i_M, j_1 \times \cdots \times j_N} &\define&
 a_{i_1, \cdots, i_M, j_1 \times \cdots \times j_N} \nonumber \\
& &+ b_{i_1, \cdots, i_M, j_1 \times \cdots \times j_N}. 
\end{eqnarray}
On the other hand, for tensors $\mathcal{A} \define (a_{i_1, \cdots, i_M, j_1, \cdots,j_N}) \in \mathbb{C}^{I_1 \times \cdots \times I_M\times
J_1 \times \cdots \times J_N}$ and $\mathcal{B} \define (b_{j_1, \cdots, j_N, k_1, \cdots,k_L}) \in \mathbb{C}^{J_1 \times \cdots \times J_N\times K_1 \times \cdots \times K_L}$, according to~\cite{liang2019further}, the \emph{Einstein product} (or simply referred to as \emph{tensor product} in this work) $\mathcal{A} \star_{N} \mathcal{B} \in  \mathbb{C}^{I_1 \times \cdots \times I_M\times
K_1 \times \cdots \times K_L}$ is given by 
\begin{eqnarray}\label{eq: Einstein product definition}
\lefteqn{(\mathcal{A} \star_{N} \mathcal{B} )_{i_1, \cdots, i_M,k_1 \times \cdots \times k_L} \define} \nonumber \\ &&\sum\limits_{j_1, \cdots, j_N} a_{i_1, \cdots, i_M, j_1, \cdots,j_N}b_{j_1, \cdots, j_N, k_1, \cdots,k_L}. 
\end{eqnarray}
Note that we will often abbreviate a tensor product $\mathcal{A} \star_{N} \mathcal{B}$ to ``$\mathcal{A} \hspace{0.05cm}\mathcal{B}$'' for notational simplicity in the rest of the paper. 
This tensor product will be reduced to the standard matrix multiplication as $L$ $=$ $M$ $=$ $N$ $=$ $1$. Other simplified situations can also be extended as tensor–vector product ($M >1$, $N=1$, and $L=0$) and tensor–matrix product ($M>1$ and $N=L=1$). In analogy to matrix analysis, we define some basic tensors and elementary tensor operations as follows. 

\begin{definition}\label{def: zero tensor}
A tensor whose entries are all zero is called a \emph{zero tensor}, denoted by $\mathcal{O}$. 
\end{definition}

\begin{definition}\label{def: identity tensor}
An \emph{identity tensor} $\mathcal{I} \in  \mathbb{C}^{I_1 \times \cdots \times I_N\times
J_1 \times \cdots \times J_N}$ is defined by 
\begin{eqnarray}\label{eq: identity tensor definition}
(\mathcal{I})_{i_1 \times \cdots \times i_N\times
j_1 \times \cdots \times j_N} \define \prod_{k = 1}^{N} \delta_{i_k, j_k},
\end{eqnarray}
where $\delta_{i_k, j_k} \define 1$ if $i_k  = j_k$; otherwise $\delta_{i_k, j_k} \define 0$.
\end{definition}
In order to define \emph{Hermitian} tensor, the \emph{conjugate transpose operation} (or \emph{Hermitian adjoint}) of a tensor is specified as follows.  
\begin{definition}\label{def: tensor conjugate transpose}
Given a tensor $\mathcal{A} \define (a_{i_1, \cdots, i_M, j_1, \cdots,j_N}) \in \mathbb{C}^{I_1 \times \cdots \times I_M\times J_1 \times \cdots \times J_N}$, its conjugate transpose, denoted by
$\mathcal{A}^{H}$, is defined by
\begin{eqnarray}\label{eq:tensor conjugate transpose definition}
(\mathcal{A}^H)_{ j_1, \cdots,j_N,i_1, \cdots, i_M}  \define  
\overline{a_{i_1, \cdots, i_M,j_1, \cdots,j_N}},
\end{eqnarray}
where the overline notion indicates the complex conjugate of the complex number $a_{i_1, \cdots, i_M,j_1, \cdots,j_N}$. If a tensor $\mathcal{A}$ satisfies $ \mathcal{A}^H = \mathcal{A}$, then $\mathcal{A}$ is a \emph{Hermitian tensor}. 
\end{definition}
\begin{definition}\label{def: unitary tensor}
Given a tensor $\mathcal{A} \define (a_{i_1, \cdots, i_M, j_1, \cdots,j_M}) \in \mathbb{C}^{I_1 \times \cdots \times I_M\times J_1 \times \cdots \times J_M}$, if
\begin{eqnarray}\label{eq:unitary tensor definition}
\mathcal{A}^H \star_M \mathcal{A} = \mathcal{A} \star_M \mathcal{A}^H = \mathcal{I} \in \mathbb{C}^{I_1 \times \cdots \times I_M\times J_1 \times \cdots \times J_M},
\end{eqnarray}
then $\mathcal{A}$ is a \emph{unitary tensor}. 
\end{definition}
\begin{definition}\label{def: inverse of a tensor}
Given a \emph{square tensor} $\mathcal{A} \define (a_{i_1, \cdots, i_M, j_1, \cdots,j_M}) \in \mathbb{C}^{I_1 \times \cdots \times I_M\times I_1 \times \cdots \times I_M}$, if there exists $\mathcal{X} \in \mathbb{C}^{I_1 \times \cdots \times I_M\times I_1 \times \cdots \times I_M}$ such that 
\begin{eqnarray}\label{eq:tensor invertible definition}
\mathcal{A} \star_M \mathcal{X} = \mathcal{X} \star_M \mathcal{A} = \mathcal{I},
\end{eqnarray}
then $\mathcal{X}$ is the \emph{inverse} of $\mathcal{A}$. We usually write $\mathcal{X} \define \mathcal{A}^{-1}$ thereby. 
\end{definition}

We also list other crucial tensor operations here. The \emph{trace} of a tensor is equivalent to the summation of all diagonal entries such that 
\begin{eqnarray}\label{eq: tensor trace def}
\mathrm{Tr}(\mathcal{A}) \define \sum\limits_{1 \leq i_j \leq I_j,\hspace{0.05cm}j \in [N]} \mathcal{A}_{i_1, \cdots, i_M,i_1, \cdots, i_M}.
\end{eqnarray}
The \emph{inner product} of two tensors $\mathcal{A}$, $\mathcal{B} \in \mathbb{C}^{I_1 \times \cdots \times I_M\times J_1 \times \cdots \times J_N}$ is given by 
\begin{eqnarray}\label{eq: tensor inner product def}
\langle \mathcal{A}, \mathcal{B} \rangle \define \mathrm{Tr}\left(\mathcal{A}^H \star_M \mathcal{B}\right).
\end{eqnarray}
According to Eq.~\eqref{eq: tensor inner product def}, the \emph{Frobenius norm} of a tensor $\mathcal{A}$ is defined by 
\begin{eqnarray}\label{eq:Frobenius norm}
\left\Vert \mathcal{A} \right\Vert_{\mathrm{F}} \define \sqrt{\langle \mathcal{A}, \mathcal{A} \rangle}.
\end{eqnarray}
In later sections, we use $\left\Vert \mathcal{A} \right\Vert_{\alpha}$ to represent the general unitary invariant norm via the gauge function $\alpha$ defined in~\cite{chang2022general}. 

Below, we will define tensor-valued processes and recall important notions about generic chaining. Let $\mathfrak{X}$ be a normed linear space made by the collection of tensors with the same dimensions and let $(T,d)$ be a metric space. We assume that the cardinality $|T|$ of $T$ is finite.  A tensor-valued process is defined as
\begin{eqnarray}\label{eq:tensor-valued process def}
\mathcal{X}_t &:& \mathbb{R} \times \Omega \rightarrow  \mathfrak{X}\nonumber \\
&& (t, \omega) \rightarrow \mathcal{X}_t (\omega).
\end{eqnarray}
The diameter of $T$ with respect to $d$ is defined as
\begin{eqnarray}\label{eq:T diameter def}
\Delta_d(T) &=& \sup\limits_{s,t \in T} d(s,t).
\end{eqnarray}
We define that an $\mathfrak{X}$-valued process $(\mathcal{X}_t)_{t \in T}$ is exponential tail with parameter $\beta$ under the metric space $(T,d)$ as
\begin{eqnarray}\label{eq:exp tail def}
\mathrm{P}(\left\Vert \mathcal{X}_t -  \mathcal{X}_s \right\Vert_\alpha \geq u d(s,t)) \leq 2 \exp(-u^\beta),
\end{eqnarray}
where $\beta \geq 0$. 

A sequence of set $\mathfrak{T} = (T_n)_{n \geq 0}$ of subsets of $T$ is called admissible if $|T_0| = 1$ and $|T_n| \leq 2^{2^n}$ for all $n \geq 1$. For any $0 < \beta < \infty$, the $\gamma_{\beta}$ functional of $(T,d)$ is defined by
\begin{eqnarray}\label{eq:gamma functional def}
\gamma_{\beta}(T,d) &=& \inf\limits_{\mathfrak{T}} \sup\limits_{t \in T}\sum\limits_{n=0}^{\infty}2^{n/\beta} d(t, T_n),
\end{eqnarray}
where the infimum is taken over all admissible sequences and $d(t, T_n)$ is obtained by $\inf\limits_{s \in T_n} d(t, s)$.  

\section{Suprema and Tail Estimates for TRP}\label{sec:Supremum and Tail Estimates for TRP} 

In this section, we will consider two types of tail bound analysis based on generic chaining. The first is to establish tail bounds for suprema of processes with an exponential tail, see Section~\ref{sec:Tail Bounds for Suprema of Processes with an Exponential Tail}. The second is to establish tail bounds for suprema of processes with mixed exponential tails, see Section~\ref{sec:Tail Bounds for Suprema of Processes with Mixed Exponential Tails}.

\subsection{Tail Bounds for Suprema of Processes with an Exponential Tail}\label{sec:Tail Bounds for Suprema of Processes with an Exponential Tail}

In order to consider $p$-th moments, we have to define the following truncated $\gamma_{\beta}$ functional. For a given $1 \leq p < \infty$, we use the symbol $n'$ to represent the integer obtained by $n' = \lfloor \log_2 (p) \rfloor$, where $\lfloor \cdot \rfloor$ represents the integer part. We define the truncated $\gamma_{\beta}$ functional with respect to the power $p$ as 
\begin{eqnarray}\label{eq:truncated gamma functional def}
\gamma_{\beta, p}(T,d) &=& \inf\limits_{\mathfrak{T}} \sup\limits_{t \in T}\sum\limits_{n=n'}^{\infty}2^{n/\beta} d(t, T_n),
\end{eqnarray} 
where $n' = \lfloor \log_2 (p) \rfloor$. By comparing definitions given by Eq.~\eqref{eq:truncated gamma functional def} and Eq.~\eqref{eq:gamma functional def}, we have $\gamma_{\beta, p}(T,d) \leq \gamma_{\beta}(T,d)$ for all $1 leq p < \infty$ and $\gamma_{\beta, 1}(T,d)  =  \gamma_{\beta}(T,d)$. The sequence $\mathfrak{T}$ that achieves the infimum is named as \emph{optimal sequence}.

We have to prepare the following three lemmas which will be used in every generic chaining argument later. 

\begin{lemma}\label{lma:A3}
Given $1 \leq p < \infty$, set $n' = \lfloor \log_2 (p) \rfloor$ and let $(\mathcal{X}_t)_{t \in T}$ be a collection of random tensors. If $|T| \leq 2^{2^{n'}}$, we have
\begin{eqnarray}\label{eq1:lma:A3}
\left(\mathbb{E} \sup\limits_{t \in T} \left\Vert \mathcal{X}_t \right\Vert_{\alpha}^{p} \right)^{1/p}
\leq 2 \sup\limits_{t \in T} \left(\mathbb{E} \left\Vert \mathcal{X}_t \right\Vert_{\alpha}^{p} \right)^{1/p}.
\end{eqnarray}
\end{lemma}
\textbf{Proof:}
Because $|T| \leq 2^{2^{n'}}$, we have 
\begin{eqnarray}
\mathbb{E}\sup\limits_{t \in T} \left\Vert \mathcal{X}_t \right\Vert_{\alpha}^{p} 
\leq \sum\limits_{t \in T} \mathbb{E} \left\Vert \mathcal{X}_t \right\Vert_{\alpha}^{p}
\leq |T| \sup\limits_{t \in T} \mathbb{E} \left\Vert \mathcal{X}_t \right\Vert_{\alpha}^{p} 
\leq 2^p \sup\limits_{t \in T} \mathbb{E} \left\Vert \mathcal{X}_t \right\Vert_{\alpha}^{p}.
\end{eqnarray}
$\hfill \Box$

\begin{lemma}\label{lma:A4}
Given $1 \leq p < \infty$, $0 < \beta \infty$, $u \geq 2^{1/\beta}$, and set $n' = \lfloor \log_2 (p) \rfloor$. For every $n > n'$, let $(E_{j,n})_{j \in J_n}$ be a collection of events satisfying 
\begin{eqnarray}\label{eq1:lma:A4}
\mathrm{P}(E_{j,n}) \leq 2 \exp(-2^n u^\beta),
\end{eqnarray}
where $j \in J_n$. If $|J_n| \leq 2^{2^{n+1}}$, we have
\begin{eqnarray}\label{eq2:lma:A4}
\mathrm{P}\left(\bigcup\limits_{n > n'}\bigcup\limits_{j \in J_n} E_{j,n}\right) \leq C_1  \exp(-p u^\beta/4).
\end{eqnarray}
\end{lemma}
\textbf{Proof:}
From union bound and the fact that $u^{\beta} \geq 2$, we have
\begin{eqnarray}\label{eq2:lma:A4}
\mathrm{P}\left(\bigcup\limits_{n > n'}\bigcup\limits_{j \in J_n} E_{j,n}\right)&\leq& \sum\limits_{n > n'}
2^{2^{n+1}}2 \exp(-2^n u^\beta) \nonumber \\
&=& 2 \sum\limits_{n > n'} \exp(2(\log 2)2^n)\exp(-2^n u^\beta) \nonumber \\
&\leq& 2 \sum\limits_{n > n'} \exp((\log 2 -1)2^n u^\beta)  \nonumber \\
&=& \exp(-2^{n'}u^{\beta}/2)\sum\limits_{n > n'} \exp((\log 2 -1)2^n u^\beta + 2^{n'}u^{\beta}/2) \nonumber \\
&\leq& \exp(-2^{n'}u^{\beta}/2)\sum\limits_{n > 0} \exp((\log 2 -1)2^n u^\beta + 2^{n}u^{\beta}/4).
\end{eqnarray}
Note that 
\begin{eqnarray}
\sum\limits_{n > 0} \exp((\log 2 -1)2^n u^\beta + 2^{n}u^{\beta}/4) 
&=& \sum\limits_{n > 0} \exp(2^{n}u^{\beta}(\log2 - 0.75)) \nonumber \\
&\leq_1&  \sum\limits_{n > 0} \exp(2^{n+1}(\log2 - 0.75)) \nonumber \\
&\leq& \sum\limits_{n > 0} \exp(2n(\log2 - 0.75)) \nonumber \\
&=& C_1.
\end{eqnarray}
where we apply $u \geq 2^{1/\beta}$ at $\leq_1$ and $(\log2 - 0.75) < 0$. Then, we have Eq.~\eqref{eq2:lma:A4} since $2^{n'} \geq p/2$. 
$\hfill \Box$

\begin{lemma}\label{lma:A5}
Given a positive random variable $x$, ratio $r \geq 0$, boundary constant $u_b >0$, $1 \leq p < \infty$ and $0 < \beta < \infty$, if we have
\begin{eqnarray}\label{eq1:lma:A5}
\mathrm{P}\left(x > r u\right) \leq d \exp(-pu^{\beta}/4),
\end{eqnarray}  
where $u \geq u_b$. Then, we have
\begin{eqnarray}\label{eq2:lma:A5}
\left(\mathbb{E}x^p\right)^{1/p} \leq r(c_{\beta}d+u_b),
\end{eqnarray}
where the constant $c_{\beta}$ depends only on $\beta$.
\end{lemma}
\textbf{Proof:}
By integration by parts and a change of variable, we have
\begin{eqnarray}\label{eq3:lma:A5}
\mathbb{E}x^p &=& \int_0^{\infty}px^{p-1} \mathrm{P}(x > u)du \nonumber \\
&=& r^p  \int_0^{\infty}pv^{p-1} \mathrm{P}(x > rv)dv \nonumber \\
&\leq_1&  r^p  \left(\int_{u_b}^{\infty}pv^{p-1} d \exp(-pu^{\beta}/4) dv + \int_0^{u_b}pv^{p-1}dv\right) \nonumber \\
&=&  r^p  \left(d\int_{u_b}^{\infty}pv^{p-1} \exp(-pu^{\beta}/4) dv + u_b^{p}\right),
\end{eqnarray}
where we apply Eq.~\eqref{eq1:lma:A5} in $\leq_1$. We need to analyze the term $\int_{u_b}^{\infty}pv^{p-1} \exp(-pu^{\beta}/4) dv$ further. Because we have 
\begin{eqnarray}
\int_{u_b}^{\infty}pv^{p-1} \exp(-pu^{\beta}/4) dv &\leq& \int_{0}^{\infty}pv^{p-1} \exp(-pu^{\beta}/4) dv
\nonumber \\
&=_1&  \frac{2p^{1 -p/\beta }2^{p/\beta}}{\beta} \int_0^{\infty}w^{2p/\beta -1}e^{-w^2/2}dw \nonumber \\
&=&  \frac{2\sqrt{2 \pi}p^{1 -p/\beta }2^{p/\beta}}{2\beta} \mathbb{E}|y|^{2p/\beta- 1}
\end{eqnarray}
where $=_1$ is obtained by change of variable and $y$ is a standard normal random variable. From the work~\cite{winkelbauer2012moments}, we have 
\begin{eqnarray}
\mathbb{E}|y|^{2p/\beta- 1} = \frac{2^{p/\beta- 1/2}}{\sqrt{\pi}}\Gamma(p/2\beta). 
\end{eqnarray}
The result follows by combining all of the above estimations. 
$\hfill \Box$

We are ready to present the following theorem about tail bounds for suprema of tensor random processes with an exponential tail.

\begin{theorem}\label{thm:3.2}
Given $0 < \beta < \infty$ and $(\mathcal{X}_t)_{t \in T}$ satisfying Eq.~\eqref{eq:exp tail def}, then there exist constants $C_\beta, D_{\beta} > 0$, which depends only on $\beta$, such that for any $t_0 \in T$ and $1 \leq p < \infty$, we have
\begin{eqnarray}\label{eq1:thm:3.2}
\left(\mathbb{E}\sup\limits_{t \in T} \left\Vert \mathcal{X}_t - \mathcal{X}_{t_0} \right\Vert_{\alpha}^p \right)^{1/p} &\leq& C_{\beta} \gamma_{\beta, p}(T,d) + 2 \sup\limits_{t \in T}\left(\mathbb{E} \left\Vert \mathcal{X}_t - \mathcal{X}_{t_0} \right\Vert_{\alpha}^p\right)^{1/p}.
\end{eqnarray}
Moreover, 
\begin{eqnarray}\label{eq2:thm:3.2}
\mathrm{P}\left(\sup\limits_{t \in T} \left\Vert \mathcal{X}_t - \mathcal{X}_{t_0} \right\Vert_{\alpha} \geq 
e^{1/\beta}\left(C_{\beta}\gamma_{\beta, p}(T,d)+uD_{\beta}\Delta_d(T)\right)\right) &\leq& \exp(-u^{\beta}/\beta).
\end{eqnarray}
\end{theorem}
\textbf{Proof:}
% \pi <-> \varpi
Let $\mathfrak{T} = (T_n)_{n \geq 0}$ be an optimal admissible sequence for $\gamma_{\beta, p}(T,d)$ and let 
$\varpi_n : T \rightarrow T_n$ be a sequence of functions defined as $\varpi_n(t) = \arg\min\limits_{s \in T_n}d(s,t)$. By setting $n' = \lfloor \log_2 (p) \rfloor$, we have the following decomposition
\begin{eqnarray}\label{eq3:thm:3.2}
\left(\mathbb{E}\sup\limits_{t \in T} \left\Vert \mathcal{X}_t - \mathcal{X}_{t_0}\right\Vert_{\alpha}^p\right)^{1/p}
&\leq& \left(\mathbb{E}\sup\limits_{t \in T} \left\Vert \mathcal{X}_t - \mathcal{X}_{\varpi_{n'}(t)}\right\Vert_{\alpha}^p\right)^{1/p} +
\left(\mathbb{E}\sup\limits_{t \in T} \left\Vert \mathcal{X}_{\varpi_{n'}(t)}- \mathcal{X}_{t_0}\right\Vert_{\alpha}^p\right)^{1/p}.
\end{eqnarray}

For the first term on R.H.S. of Eq.~\eqref{eq3:thm:3.2}, we will apply telescoping sum representation for $\mathcal{X}_t - \mathcal{X}_{\varpi_{n'}(t)}$. It is 
\begin{eqnarray}\label{eq4:thm:3.2}
\mathcal{X}_t - \mathcal{X}_{\varpi_{n'}(t)} &=& \sum\limits_{n > n'} \mathcal{X}_{\varpi_{n}(t)} -  \mathcal{X}_{\varpi_{n-1}(t)}.
\end{eqnarray}
Due to that the tensor random process $(\mathcal{X}_t)_{t \in T}$ satisfies Eq.~\eqref{eq:exp tail def}, we have 
\begin{eqnarray}\label{eq5:thm:3.2}
\mathrm{P}\left(\left\Vert \mathcal{X}_{\varpi_{n}(t)} -  \mathcal{X}_{\varpi_{n-1}(t)} \right\Vert_{\alpha}\geq  u2^{n/\beta} d(\varpi_n(t),\varpi_{n-1}(t))\right) \leq 2 \exp(-u^\beta 2^n),
\end{eqnarray}
where $n > n'$. If $u \geq 2^{1/\beta}$, we use $E_{u,p}$ to represent the following event:
\begin{eqnarray}\label{eq6:thm:3.2}
\left\Vert \mathcal{X}_{\varpi_{n}(t)} -  \mathcal{X}_{\varpi_{n-1}(t)} \right\Vert_{\alpha} \leq  u2^{n/\beta} d(\varpi_n(t),\varpi_{n-1}(t)),
\end{eqnarray}
where $\forall n > n'$ and $\forall t \in T$. Then, by Lemma~\ref{lma:A4}, we have the following probability bound for the complement of the event $E_{u,p}$:
\begin{eqnarray}\label{eq7:thm:3.2}
\mathrm{P}(E^c_{u,p}) \leq C_1 \exp(-pu^{\beta}/4).
\end{eqnarray}
If the event $E_{u,p}$ happens, we have
\begin{eqnarray}\label{eq8:thm:3.2}
\left\Vert \sum\limits_{n > n'} \mathcal{X}_{\varpi_{n}(t)} -  \mathcal{X}_{\varpi_{n-1}(t)}\right\Vert_{\alpha}
&\leq& \sum\limits_{n > n'}  \left\Vert \mathcal{X}_{\varpi_{n}(t)} -  \mathcal{X}_{\varpi_{n-1}(t)}\right\Vert_{\alpha} \nonumber \\
&\leq& u  \sum\limits_{n > n'} 2^{n/\beta} d(\varpi_n(t),\varpi_{n-1}(t))   \nonumber \\
&\leq& u(1+2^{1/\beta})\gamma_{\beta, p}(T,d).
\end{eqnarray}
Therefore, we have 
\begin{eqnarray}\label{eq9:thm:3.2}
\mathrm{P}\left(\sup\limits_{t \in T} \left\Vert \mathcal{X}_t - \mathcal{X}_{\varpi_{n'}(t)} \right\Vert_\alpha u(1+2^{1/\beta})\gamma_{\beta, p}(T,d)\right) \leq  C_1 \exp(-pu^{\beta}/4).
\end{eqnarray}
By Lemma~\ref{lma:A5}, we have
\begin{eqnarray}\label{eq10:thm:3.2}
\left(\mathbb{E}\sup\limits_{t \in T} \left\Vert \mathcal{X}_t - \mathcal{X}_{\varpi_{n'}(t)} \right\Vert_\alpha^p \right)^{1/p} \leq C_\beta \gamma_{\beta, p}(T,d).
\end{eqnarray}

For the second term on R.H.S. of Eq.~\eqref{eq3:thm:3.2}, we can estimate this by Lemma~\ref{lma:A3}. It is 
\begin{eqnarray}\label{eq11:thm:3.2}
\left(\mathbb{E}\sup\limits_{t \in T} \left\Vert \mathcal{X}_{\varpi_{n'}(t)}- \mathcal{X}_{t_0}\right\Vert_{\alpha}^p\right)^{1/p} &\leq& 2 \sup\limits_{t \in T}\left(\mathbb{E} \left\Vert \mathcal{X}_{\varpi_{n'}(t)}- \mathcal{X}_{t_0}\right\Vert_{\alpha}^p \right)^{1/p} \nonumber \\
&\leq&  2 \sup\limits_{t \in T}\left(\mathbb{E} \left\Vert \mathcal{X}_{t}- \mathcal{X}_{t_0}\right\Vert_{\alpha}^p \right)^{1/p}. 
\end{eqnarray}

Then, Eq.~\eqref{eq1:thm:3.2} follows by combining Eq.~\eqref{eq3:thm:3.2}, Eq.~\eqref{eq10:thm:3.2} and Eq.~\eqref{eq11:thm:3.2}.

Before proving Eq.~\eqref{eq2:thm:3.2}, we have to state two facts about a random variable. Their proof can be found at Chapter 7 in~\cite{foucart2013invitation}. If $X$ is a random variable satisfying 
\begin{eqnarray}
(\mathbb{E}|X|^p)^{1/p} \leq a p^{1/\beta} + b,
\end{eqnarray}
where $0<a,b,\beta < \infty$ and $p \geq 1$, then, we have
\begin{eqnarray}\label{eq12:thm:3.2}
\mathrm{P}\left(|X| \geq e^{1/\beta}(au + b)\right) \leq \exp(-u^{\beta}/\beta),
\end{eqnarray}
where $u \geq 1$. The second fact is about how to bound the random variable moment from its tail bound. If a random variable $X$ satisfies the following:
\begin{eqnarray}
\mathrm{P}(|X| \geq e^{1/\beta}au) \leq be^{-u^{\beta}/\beta},
\end{eqnarray} 
where $0 < a,b, \beta, u < \infty$; we have 
\begin{eqnarray}\label{eq13:thm:3.2}
(\mathbb{E}|X|^p)^{1/p} \leq e^{1/2e}ap^{1/\beta}\left(\sqrt{\frac{2\pi}{\beta}}b e^{\beta/12}\right)^{1/p}.
\end{eqnarray}

From the condition provided by Eq.~\eqref{eq:exp tail def} and the fact given by Eq.~\eqref{eq13:thm:3.2}, we have
\begin{eqnarray}\label{eq14:thm:3.2}
\sup\limits_{t \in T}\left(\mathbb{E} \left\Vert \mathcal{X}_t - \mathcal{X}_{t_0}\right\Vert_{\alpha}\right)^{1/p}
\leq D_{\beta}\Delta_d(T)p^{1/\beta}.
\end{eqnarray}
Then Eq.~\eqref{eq2:thm:3.2} is obtained by applying Eq.~\eqref{eq14:thm:3.2} to Eq.~\eqref{eq1:thm:3.2} with Eq.~\eqref{eq12:thm:3.2}.
$\hfill \Box$

We use the symbol $\preceq$ in $\mathcal{X}  \preceq  \mathcal{Y}$ to represent that the tensor $\mathcal{Y}- \mathcal{X}$ is a positive definite tensor. 

Since Theorem~\ref{thm:3.2} does not need any independence assumptions on the increments of the tensor process $(\mathcal{X}_t)_{t \in T}$, we can apply Theorem~\ref{thm:3.2} to martingale context. If $\left(\mathcal{X}_i\right)_{i=0,1,2,\cdots,n} \in \mathbb{C}^{I_1 \times \cdots \times I_M  \times I_1 \times \cdots \times I_M }$ is a discrete-time Hermitian tensor-valued martingale and its $i$-th difference tensor is  by $\mathcal{D}_i = \mathcal{X}_i - \mathcal{X}_{i-1}$, then, from Theorem 14 (Tensor Azuma Inequality) in~\cite{chang2021convenient}, we have
\begin{eqnarray}\label{eq:Azuma inequality}
\mathrm{P}\left(\lambda_{\max}\left(\mathcal{X}_n - \mathcal{X}_{0}\right) \geq u \right) \leq 
\left(\prod\limits_{i=1}^M I_i\right)\exp\left(-\frac{u^2}{8 \sigma^2} \right),
\end{eqnarray}
where $\sigma^2 = \left\Vert \sum\limits_{i=1}^n \mathcal{D}_i^2 \right\Vert_{\mbox{\tiny max spec}}$ and $\left\Vert\cdot \right\Vert_{\mbox{\tiny max spec}}$ represents the maximum spectral norm. Combing this Tensor Azuma Inequality given by Eq.~\eqref{eq:Azuma inequality} with Theorem~\ref{thm:3.2}, we have the following corollary.

\begin{corollary}\label{cor:3.4}
Let $\mathcal{X}_t = \left(\mathcal{X}_{t,i}\right)_{i=0,1,\cdots,n}, t \in T$ be a family of discrete-time martingales
with respect to the same filtration. We define the following metric:
\begin{eqnarray}\label{eq1:cor:3.4}
d_{s,t} = \left(\left\Vert \sum\limits_{i=1}^n \mathcal{D}_i^2 \right\Vert_{\mbox{\tiny max spec}} \right)^{1/2},
\end{eqnarray}
where $\mathcal{D}_i = (\mathcal{X}_{t,i} - \mathcal{X}_{s,i}) - (\mathcal{X}_{t,i-1} - \mathcal{X}_{s,i-1})$. For any $u \geq 1$, we have
\begin{eqnarray}\label{eq2:cor:3.4}
\mathrm{P}\left(\sup\limits_{t \in T} \lambda_{\max}\left(\mathcal{X}_{t,n}- \mathcal{X}_{t,0}\right)\geq e^{1/2}\left(C_2 \gamma_{2}(T,d) + D_2 \Delta_d(T)u\right) \right) \leq e^{-u^2/2}.
\end{eqnarray}
\end{corollary}

\subsection{Tail Bounds for Suprema of Processes with Mixed Exponential Tails}\label{sec:Tail Bounds for Suprema of Processes with Mixed Exponential Tails}

The purpose of this section is to establish tail bouns for suprema of tensor processes with mixed exponential tails. We have to prepare following two lemmas which will be used in proving tail bouns for suprema of tensor processes.

\begin{lemma}\label{lma:A1}
If $X$ is a random variable satisfying
\begin{eqnarray}\label{eq1:lma:A1}
\left(\mathbb{E}|X|^p\right)^{1/p} \leq \sum\limits_{n=1}^{m}a_n p^{1/n} + a_{n+1},
\end{eqnarray} 
where $p \geq 1$ and $0 \leq a_n < \infty$; then we have 
\begin{eqnarray}\label{eq2:lma:A1}
\mathrm{P}\left(|X| \leq e\left(\sum\limits_{n=1}^{m}a_n p^{1/n} + a_{n+1}\right)\right) \leq \exp(-u),
\end{eqnarray} 
where $u \geq 1$. 
\end{lemma}
\textbf{Proof:}
The proof is a straightforward consequence of Markov’s inequality. 
$\hfill \Box$

\begin{lemma}\label{lma:A2}
If a random variable $X$ satisfies
\begin{eqnarray}\label{eq1:lma:A2}
\mathrm{P}\left(|X| \geq \sum\limits_{n=1}^{m}a_n p^{1/n} \right) \leq \exp(-u),
\end{eqnarray}
where $u \geq 0$ and $0 \leq a_n < \infty$; then, 
\begin{eqnarray}\label{eq2:lma:A2}
\left(\mathbb{E}|X|^p\right)^{1/p} \leq \sum\limits_{n=1}^{m} m a_n f_n(p) p^{1/n},
\end{eqnarray}
where $f_n(p)$ are positive real-valued functions related to the parameter $p$ and $p \geq 1$.
\end{lemma}
\textbf{Proof:}
From Eq.~\eqref{eq2:lma:A2}, we have
\begin{eqnarray}\label{eq3:lma:A2}
\mathrm{P}\left(\frac{1}{m}|X| \geq u \right) \leq \begin{cases}
e^{-u^m/a^m_m}, & \mbox{if $0  \leq u \leq a_m^m/a^{m-1}_{m-1}$};\\
~~~~\vdots~~~~, & ~~~~ \vdots ~~~~;\\
e^{-u^3/a^3_3}, & \mbox{if $a_4^4/a^3_3  \leq u \leq a_3^3/a^2_2 $};\\
e^{-u^2/a^2_2}, & \mbox{if $a_3^3/a^2_2  \leq u \leq a_2^2/a_1$};\\
e^{-u/a_1}, & \mbox{if $u \geq a_2^2/a_1$}.
\end{cases}
\end{eqnarray}
Using integration by parts and a change of variable, we have
\begin{eqnarray}\label{eq4:lma:A2}
m^{-p}\mathbb{E}|X|^p &=& p \int_0^{\infty}u^{p-1}\mathrm{P}\left(\frac{1}{m}|X| \geq u \right) du \nonumber \\
&\leq& p  \int_0^{a_m^m/a^{m-1}_{m-1}}u^{p-1}e^{-u^m/a^m_m}du + 
p \sum\limits_{n=2}^{m-1}  \int_{a_{m-n+2}^{m-n+2}/a_{m-n+1}^{m-n+1}}^{a_{m-n+1}^{m-n+1}/a^{m-n}_{m-n}}u^{p-1}e^{-u^{m-n+1}/a^{m-n+1}_{m-n+1}}du \nonumber \\
&& + p  \int_{a_2^2/a_1}^{\infty}u^{p-1}e^{-u/a_1}du \nonumber \\
&\leq& \frac{1}{m}p a_m^p \Gamma(p/m) + \sum\limits_{n=2}^{m-1}\frac{1}{n}p a_n^p \Gamma(p/n) + pa_1^p \Gamma(p)  \nonumber \\
&=&  \sum\limits_{n=1}^{m}\frac{1}{n}p a_n^p \Gamma(p/n) \nonumber  \\
&\leq_1& \sum\limits_{n=1}^{m} m a_n f_n(p) p^{1/n},
\end{eqnarray}
where $\Gamma(\cdot)$ is the gamma function and $\leq_1$ comes from estimation for Stirling’s formula. The Stirling’s formula of $\Gamma(p)$ can be expressed as  
\begin{eqnarray}\label{eq5:lma:A2}
\Gamma(p) &=& \sqrt{2\pi}p^{p - 1/2}e^{-p}e^{\mu(p)/12p},
\end{eqnarray}
where $0 \leq \mu(p) \leq 1$. This Lemma is proved by finding $f_n(p)$ to bound $\frac{1}{n}p \Gamma(p/n) $ for each $n$. For example, if $m=2$, we have 
\begin{eqnarray}
f_1(p) &=& \sqrt{2\pi p}p^p e^{-p+1/12p}, \nonumber \\
f_2(p) &=& \sqrt{\pi}e^{1/6p}(2e)^{-p/2}e^{p/2e}p^{p/2}.
\end{eqnarray}
$\hfill \Box$

With these two lemmas, we are ready to present the following theorem about tail bouns for suprema of tensor processes with mixed exponential tails.
\begin{theorem}\label{thm:3.5}
If $\left(\mathcal{X}_t\right)_{t \in T}$ has a mixed tail which satisfies the following:
\begin{eqnarray}\label{eq1:thm:3.5}
\mathrm{P}\left(\left\Vert \mathcal{X}_t - \mathcal{X}_s \right\Vert_{\alpha} \geq 
\sum\limits_{n=1}^m u^{1/n}d_n(s,t)\right) \leq 2e^{-u},
\end{eqnarray}
where $d_n(s,t)$ are metrics defined on $T$ and $u \geq 0$. Then, there is a constant $C_2 > 0$ such that for any $1 \leq p < \infty$, we have
\begin{eqnarray}\label{eq2:thm:3.5}
\left(\mathbb{E}\sup\limits_{t \in T} \left\Vert \mathcal{X}_t - \mathcal{X}_{t_0} \right\Vert_{\alpha}\right)^{1/p}
&\leq&C_2\sum\limits_{n=1}^m \gamma_n(T,d_n) + 2 \sup\limits_{t \in T}\left(\mathbb{E}\left\Vert \mathcal{X}_t - \mathcal{X}_{t_0}\right\Vert^p_{\alpha}\right)^{1/p}.
\end{eqnarray}
Accordingly, there are constants $C_2, C_3$ such that for any $u \geq 1$, we have
\begin{eqnarray}\label{eq3:thm:3.5}
\mathrm{P}\left(\sup\limits_{t \in T}\left\Vert \mathcal{X}_t - \mathcal{X}_{t_0} \right\Vert_{\alpha}\geq 
C_2 \sum\limits_{n=1}^m \gamma_n(T,d_n)  + C_3 \sum\limits_{n=1}^m u^{1/n}\Delta_{d_n}(T) \right)\leq e^{-u}. 
\end{eqnarray}
\end{theorem}
\textbf{Proof:}
We select $m$ admissible sequences of paritions $\mathfrak{Q}^{(n)} = (\bm{Q}^{(n)}_i)_{i \geq 0}$ for $n=1,2,\cdots,m$, such that 
\begin{eqnarray}\label{eq4:thm:3.5}
\sup\limits_{t \in T}\sum\limits_{i \geq 0}2^{i/n}\Delta_{d_n}(\bm{Q}^{(n)}_i(t)) \leq 2 \gamma'_n(T,d_n),
\end{eqnarray}
where $\gamma'_n(T,d_n)$ is defined as 
\begin{eqnarray}\label{eq5:thm:3.5}
\gamma'_n(T,d_n) &\define& \inf\limits_{\mathfrak{Q}^{(n)}}\sup\limits_{t \in T} \sum\limits_{i=0}^{\infty}2^{i/n} \Delta_{d_n}(\bm{Q}^{(n)}_i(t)).
\end{eqnarray}

For $i \geq 0$, let $\mathfrak{A} = (A_i)_{i \geq 0}$ be the partition generated by $\bm{Q}^{(n)}_{i-j}$ by
\begin{eqnarray}\label{eq6:thm:3.5}
A_i = \left\{\bigcap\limits_{n=1}^m Q^{(n)}:  Q^{(n)} \in \bm{Q}^{(n)}_{i-j} \right\}
\end{eqnarray}
where $j =\lceil \log_2 m\rceil$ and $\bm{Q}^{(n)}_{-j} \define \bm{Q}^{(n)}_{0}$ so that $A_0 = T$. Then, $(A_i)_{i \geq 0}$ is increasing and 
\begin{eqnarray}\label{eq7:thm:3.5}
|A_i| \leq \prod\limits_{n=1}^m |\bm{Q}^{(n)}_{i-j}| \leq 2^{2^{i-j}} \cdot 2^{2^{i-j}}\cdots 2^{2^{i-j}}  = 2^{2^{i}},
\end{eqnarray}
which shows that $\mathfrak{A}$ is admissible. For every $i \geq 0$, we define $T_i$ of $T$ by selecting exactly one point from each $A_i$. We define $\varpi_i(t): T \rightarrow T_i$ be the unique element of $T_i \bigcap A_i(t)$, then we have a sequence $(\varpi_i)_{i \geq 0}$ of maps.

By setting $i' = \lfloor \log_2 (p) \rfloor$, we have the following decomposition
\begin{eqnarray}\label{eq8:thm:3.5}
\left(\mathbb{E}\sup\limits_{t \in T} \left\Vert \mathcal{X}_t - \mathcal{X}_{t_0}\right\Vert_{\alpha}^p\right)^{1/p}
&\leq& \left(\mathbb{E}\sup\limits_{t \in T} \left\Vert \mathcal{X}_t - \mathcal{X}_{\varpi_{i'}(t)}\right\Vert_{\alpha}^p\right)^{1/p} +
\left(\mathbb{E}\sup\limits_{t \in T} \left\Vert \mathcal{X}_{\varpi_{i'}(t)}- \mathcal{X}_{t_0}\right\Vert_{\alpha}^p\right)^{1/p}.
\end{eqnarray}

We will focus the estimation for the first term on R.H.S. of Eq.~\eqref{eq8:thm:3.5} since the second term term on R.H.S. of Eq.~\eqref{eq8:thm:3.5} can be derived similarly in the proof of Theorem~\ref{thm:3.2}. For the frst term, we can apply the telescoping sum as 
\begin{eqnarray}
\mathcal{X}_t - \mathcal{X}_{\varpi_{i'}(t)} = \sum\limits_{i > i'} \mathcal{X}_{\varpi_i(t)} -  \mathcal{X}_{\varpi_{i-1}(t)}.
\end{eqnarray}
From condition of mixed tail given by Eq.~\eqref{eq1:thm:3.5}, for $i > i'$ and $u \geq 0$, we have
\begin{eqnarray}
\mathrm{P}\left(\left\Vert \mathcal{X}_{\varpi_i(t)} - \mathcal{X}_{\varpi_{i-1}(t)} \right\Vert_{\alpha} \geq 
\sum\limits_{n=1}^m u^{1/n}2^{i/n}d_n(\varpi_i(t),\varpi_{i-1}(t))\right) \leq 2e^{-u2^i}.
\end{eqnarray}
Let $E_{u,p}$ represent the event 
\begin{eqnarray}
\left\Vert \mathcal{X}_{\varpi_i(t)} - \mathcal{X}_{\varpi_{i-1}(t)} \right\Vert_{\alpha} \geq 
\sum\limits_{n=1}^m u^{1/n}2^{i/n}d_n(\varpi_i(t),\varpi_{i-1}(t)), 
\end{eqnarray}
where $i > i'$ and $t \in T$. If the event $E_{u,p}$ happens, we have
\begin{eqnarray}
\left\Vert \sum\limits_{i > i'} \mathcal{X}_{\varpi_i(t)} - \mathcal{X}_{\varpi_{i-1}(t)}\right\Vert_{\alpha}
&\leq &\sum\limits_{i > i'} \left\Vert \mathcal{X}_{\varpi_i(t)} - \mathcal{X}_{\varpi_{i-1}(t)} \right\Vert_{\alpha} \nonumber  \\
&\leq& \sum\limits_{n=1}^m u^{1/n}\sum\limits_{i > i'}2^{i/n}d_n(\varpi_i(t),\varpi_{i-1}(t)).
\end{eqnarray}

By our choice of $\mathfrak{Q}^{(n)}$, we can find constrants $h_n$ for $n=1,2,\cdots,m$ such that 
\begin{eqnarray}
\sum\limits_{i > i'}2^{i/n}d_n(\varpi_i(t),\varpi_{i-1}(t)) \leq h_n \gamma'_n (T,d_n).
\end{eqnarray}
Due to the fact that $\sup\limits_{t \in T}\left\Vert \mathcal{X}_{t} - \mathcal{X}_{\varpi_{n'}(t)} \right\Vert_{\alpha} \leq \sum\limits_{n=1}^m h_n \gamma'_n (T,d_n)$, we have
\begin{eqnarray}
\mathrm{P}\left( \sup\limits_{t \in T}\left\Vert \mathcal{X}_{t} - \mathcal{X}_{\varpi_{n'}(t)} \right\Vert_{\alpha}>  C'_2 u \sum\limits_{n=1}^m \gamma'_n (T,d_n) \right) \leq c \exp(-pu/4),
\end{eqnarray}
where $C'_2=\max\limits_{n=1,2,\cdots,m}h_n$. By Lemma~\ref{lma:A5} and the fact $\gamma_n (T,d_n) \leq \gamma'_n (T,d_n) \leq  g(n) \gamma_n (T,d_n)$, where $g(n)$ is a constant depending only on $n$~\cite{}; we have
\begin{eqnarray}
\left(\mathbb{E} \sup\limits_{t \in T} \left\Vert \mathcal{X}_{t} - \mathcal{X}_{\varpi_{n'}(t)} \right\Vert_{\alpha}^p\right)^{1/p} \leq C'_2 u \sum\limits_{n=1}^m \gamma'_n (T,d_n) \leq C_2 u \sum\limits_{n=1}^m \gamma_n (T,d_n).
\end{eqnarray}
Then, we have establish Eq.~\eqref{eq2:thm:3.5}.

For the tail bound, Lemma~\ref{lma:A2} and the mixed tail condition provided by Eq.~\eqref{eq1:thm:3.5} give us 
\begin{eqnarray}\label{eq9:thm:3.5}
\sup\limits_{t \in T}\left(\mathbb{E}\left\Vert \mathcal{X}_t - \mathcal{X}_{t_0}\right\Vert_{\alpha}^p \right)^{1/p}
\leq C_3 \sum\limits_{n=1}^m p^{1/n}\Delta_{d_n} (T).
\end{eqnarray}
Eq.~\eqref{eq3:thm:3.5} can be obtained by applying Eq.~\eqref{eq9:thm:3.5} to Eq.~\eqref{eq2:thm:3.5} with Lemma~\ref{lma:A1}.
$\hfill \Box$

\section{Application: Compressed Sensing}\label{sec:Application: Compressed Sensing} 

In this section, we will apply Theorem~\ref{thm:3.2} to high-dimensional compressed sensing. We begin with some basic notations for high-dimensional compressed sensing based on tensors. Note that $\left\Vert \cdot \right\Vert_{\mathrm{F}}$ is the standard Frobenius norm.

We say that a tensor $\mathcal{X} \in \mathbb{C}^{J_1 \times \cdots \times J_N}$ is $\xi$-sparse if 
\begin{eqnarray}
\left\Vert \mathcal{X} \right\Vert_0  \define |\{(j_1, j_2, \cdots, j_N): x_{j_1, j_2, \cdots, j_N} \neq 0 \}| \leq \xi.
\end{eqnarray}
For a given $\xi \in \mathbb{N}$, the \emph{$\xi$-th restricted isometry constant} $\tau_{\xi}$ of an $I_1 \times \cdots \times I_M \times J_1 \times \cdots \times J_N$ tensor $\mathcal{A}$ with the restrction $\prod\limits_{i=1}^M I_i \leq \prod\limits_{j=1}^N J_j$ is the smallest constant $\tau$ such that 
\begin{eqnarray}
(1 - \tau)\left\Vert \mathcal{X} \right\Vert_{\mathrm{F}}^2 \leq \left\Vert \mathcal{A} \star_N \mathcal{X} \right\Vert_{\mathrm{F}}^2  \leq (1 + \tau)\left\Vert \mathcal{X} \right\Vert_{\mathrm{F}}^2,
\end{eqnarray}
for all $\xi$-sparse tensor $\mathcal{X} \in \mathbb{C}^{J_1 \times \cdots \times J_N}$. We also can define a set for $\xi$-sparse tensors $\mathcal{X}$ as 
\begin{eqnarray}
S_{\xi;(J_1,\cdots,J_N)} &=& \left\{\mathcal{X} \in \mathbb{C}^{J_1 \times \cdots \times J_N}: \left\Vert \mathcal{X} \right\Vert_{\mathrm{F}} = 1, \left\Vert \mathcal{X} \right\Vert_0  \leq \xi \right\},
\end{eqnarray}
then, we have
\begin{eqnarray}
\tau_{\xi}(\mathcal{A}) = \sup\limits_{\mathcal{X} \in S_{\xi;(J_1,\cdots,J_N)}} \left\vert \left\Vert \mathcal{A} \star_N \mathcal{X} \right\Vert_{\mathrm{F}}^2 -1 \right\vert.
\end{eqnarray}
The restricted isometry constant $\tau_{\xi}(\mathcal{A}) $ is an important parameter in compressed sensing.

Let $[u_{i_1,\cdots,i_N,j_1,\cdots,j_N}] \in \mathcal{U} \in J_1 \times \cdots \times J_N \times J_1 \times \cdots \times J_N$ be a unitary tensor and we have some positive contant $\Upsilon \geq 1$ such that 
\begin{eqnarray}
\sup\limits_{i_1,\cdots,i_n,j_1,\cdots,j_n} \sqrt{\prod\limits_{j=1}^N J_j}|u_{i_1,\cdots,i_n,j_1,\cdots,j_n}| \leq \Upsilon.
\end{eqnarray}
We consider a sequence $(\theta_{j_1,\cdots,j_N})_{1 \leq j_i \leq J_j}$ for $i=1,2,\cdots,N$ of i.i.d. copies of the random selector, i.e., a Bernoulli random variable that satisfies 
\begin{eqnarray}
\mathrm{P}\left(\theta = 1\right) = \frac{\prod\limits_{i=1}^M I_i}{\prod\limits_{j=1}^N J_j}.
\end{eqnarray}
Let $\Pi = \{(j_1,\cdots,j_N):\theta_{j_1,\cdots,j_N} = 1\}$ be the random set of selected indices with
expected cardinality as $\mathbb{E}|\Pi| = \prod\limits_{i=1}^M I_i$. Let $\mathcal{U}_{\Pi}$ be a sampled tensor obtained by the unitary tensor $\mathcal{U}$ as
\begin{eqnarray}
\mathcal{U}_{\Pi} &=& \sqrt{\frac{\prod\limits_{j=1}^N J_j}{\prod\limits_{i=1}^M I_i}} \mathcal{S}_{\Pi} \star_N \mathcal{U},
\end{eqnarray}
where $\mathcal{S}_{\Pi}$ is the operator which restricts a tensor to its entries in $\Pi$ of $\mathcal{U}$. 

We are ready to present following theorem about the restricted isometry constant $\tau_{\xi}$.
\begin{theorem}\label{thm:4.1}
If there are universal constants $C_4, C_5 > 0$ such that for any given $\xi \in \mathbb{N}$ and any $0< \tau, \eta < 1$ satisfying  
\begin{eqnarray}\label{eq1:thm:4.1}
\prod\limits_{i=1}^M I_i \geq \xi \Upsilon^2 \tau^2 \max\{C_4 \log^2\xi \log\left(\prod\limits_{i=1}^M I_i \right) \log\left(\prod\limits_{j=1}^N J_j\right), C_5 \log(\eta^{-1})\},
\end{eqnarray}
then, 
\begin{eqnarray}\label{eq2:thm:4.1}
\mathrm{P}\left(\tau_{\xi}(\mathcal{U}_{\Pi}) \geq \tau\right) \leq \eta.
\end{eqnarray} 
\end{theorem}
\textbf{Proof:}
Let $\mathcal{U}_{j_1,\cdots,j_N}$ be the $(j_1,\cdots,j_N)-$ th projection tensor of $\mathcal{U}$. For every    
tensor $\mathcal{X} \in S_{\xi;(J_1,\cdots,J_N)}$, we define following function
\begin{eqnarray}\label{eq3:thm:4.1}
f_{\mathcal{X}}(\theta_{j_1,\cdots,j_N}) = \theta_{j_1,\cdots,j_N}\sqrt{\frac{\prod\limits_{j=1}^N J_j}{\prod\limits_{i=1}^M I_i}}\langle \mathcal{U}_{j_1,\cdots,j_N}, \mathcal{X} \rangle.
\end{eqnarray}
Because $\mathcal{U}$ is an unitary tensor, we have
\begin{eqnarray}
\sum\limits_{j_1=1,\cdots,j_N=1}^{J_1,\cdots,J_N}\mathbb{E}f_{\mathcal{X}}(\theta_{j_1,\cdots,j_N}) 
&=& \frac{\prod\limits_{j=1}^N J_j}{\prod\limits_{i=1}^M I_i}\sum\limits_{j_1=1,\cdots,j_N=1}^{J_1,\cdots,J_N}\mathbb{E}(\theta_{j_1,\cdots,j_N})\left\vert \langle \mathcal{U}_{j_1,\cdots,j_N}, \mathcal{X}\rangle\right\vert^2 \nonumber \\
&=&  \sum\limits_{j_1=1,\cdots,j_N=1}^{J_1,\cdots,J_N}\left\vert \langle \mathcal{U}_{j_1,\cdots,j_N}, \mathcal{X}\rangle\right\vert^2 \nonumber \\
&=& \left\Vert \mathcal{U}\star_N \mathcal{X}\right\Vert_2^2 = 1.
\end{eqnarray} 
We can express $\tau_{\xi}(\mathcal{U}_{\Pi})$ as 
\begin{eqnarray}
\tau_{\xi}(\mathcal{U}_{\Pi}) &=& \sup\limits_{\mathcal{X} \in S_{\xi;(J_1,\cdots,J_N)}}\left\vert \sum\limits_{j_1=1,\cdots,j_N=1}^{J_1,\cdots,J_N}\left(f^2_{\mathcal{X}}(\theta_{j_1,\cdots,j_N})- \mathbb{E}f^2_{\mathcal{X}}(\theta_{j_1,\cdots,j_N})\right)\right\vert.
\end{eqnarray}

Let $\epsilon_{j_1,\cdots,j_N}$ be a Rademacher sequence, i.e., a sequence of independent
symmetric Bernoulli random variables, and $1 \leq p < \infty$, we have
\begin{eqnarray}\label{eq:4.5}
\left(\mathbb{E}\tau^p_{\xi}(\mathcal{U}_{\Pi})\right)^{1/p} \leq 2 \left(\mathbb{E} \mathbb{E}_{\epsilon}\sup\limits_{\mathcal{X}\in S_{\xi;(J_1,\cdots,J_N)}} \left\vert \sum\limits_{j_1=1,\cdots,j_N=1}^{J_1,\cdots,J_N}\epsilon_{j_1,\cdots,j_N}f_{\mathcal{X}}(\theta_{j_1,\cdots,j_N})\right\vert^p \right)^{1/p}.
\end{eqnarray}
By Hoeffding's inequality, we have
\begin{eqnarray}
\mathrm{P}\left(\sum\limits_{j_1=1,\cdots,j_N=1}^{J_1,\cdots,J_N}\epsilon_{j_1,\cdots,j_N}\left(
(f^2_{\mathcal{X}}(\theta_{j_1,\cdots,j_N}) - f^2_{\mathcal{Y}}(\theta_{j_1,\cdots,j_N}))^2\right)^{1/2}\right) \leq \exp(-u^2/2).  
\end{eqnarray}
Because we have
\begin{eqnarray}
\lefteqn{\left(\sum\limits_{j_1=1,\cdots,j_N=1}^{J_1,\cdots,J_N}\left(
f^2_{\mathcal{X}}(\theta_{j_1,\cdots,j_N}) - f^2_{\mathcal{Y}}(\theta_{j_1,\cdots,j_N})\right)^2\right)^{1/2}
}\nonumber \\
&=&\left(\sum\limits_{j_1=1,\cdots,j_N=1}^{J_1,\cdots,J_N}\left(
f_{\mathcal{X}}(\theta_{j_1,\cdots,j_N}) - f_{\mathcal{Y}}(\theta_{j_1,\cdots,j_N})\right)^2\left(
f_{\mathcal{X}}(\theta_{j_1,\cdots,j_N}) + f_{\mathcal{Y}}(\theta_{j_1,\cdots,j_N})\right)^2\right)^{1/2}  \nonumber \\
&\leq&  2 \sup\limits_{\mathcal{Z} \in  S_{\xi;(J_1,\cdots,J_N)}}\left(\sum\limits_{j_1=1,\cdots,j_N=1}^{J_1,\cdots,J_N}f_{\mathcal{Z}}(\theta_{j_1,\cdots,j_N})\right)^{1/2} \max\limits_{\substack{1\leq j_i \leq J_i,\\ \mbox{\tiny for $1\leq i \leq N$}}}\left\vert f_{\mathcal{X}}(\theta_{j_1,\cdots,j_N}) - f_{\mathcal{Y}}(\theta_{j_1,\cdots,j_N})\right\vert,
\end{eqnarray}
then, we can hve the following subgaussian process
\begin{eqnarray}
\mathcal{X} \rightarrow \sum\limits_{j_1=1,\cdots,j_N=1}^{J_1,\cdots,J_N}\epsilon_{j_1,\cdots,j_N}
 f^2_{\mathcal{X}}(\theta_{j_1,\cdots,j_N}),
\end{eqnarray}
with respect to the metric
\begin{eqnarray}
\lefteqn{d(\mathcal{X}, \mathcal{Y})}
\nonumber \\
&=& 2 \sup\limits_{\mathcal{Z} \in  S_{\xi;(J_1,\cdots,J_N)}}\left(\sum\limits_{j_1=1,\cdots,j_N=1}^{J_1,\cdots,J_N}f_{\mathcal{Z}}(\theta_{j_1,\cdots,j_N})\right)^{1/2} \max\limits_{\substack{1\leq j_i \leq J_i,\\ \mbox{\tiny for $1\leq i \leq N$}}}\left\vert f_{\mathcal{X}}(\theta_{j_1,\cdots,j_N}) - f_{\mathcal{Y}}(\theta_{j_1,\cdots,j_N})\right\vert \nonumber \\
&=&  2 \sup\limits_{\mathcal{Z} \in  S_{\xi;(J_1,\cdots,J_N)}}\left(\sum\limits_{j_1=1,\cdots,j_N=1}^{J_1,\cdots,J_N}f_{\mathcal{Z}}(\theta_{j_1,\cdots,j_N})\right)^{1/2}d_{\theta}(\mathcal{X}, \mathcal{Y}),
\end{eqnarray}
where we represent the metric $d_{\theta}(\mathcal{X}, \mathcal{Y})$ as 
\begin{eqnarray}
d_{\theta}(\mathcal{X}, \mathcal{Y})&=&\max\limits_{\substack{1\leq j_i \leq J_i,\\ \mbox{\tiny for $1\leq i \leq N$}}}\left\vert f_{\mathcal{X}}(\theta_{j_1,\cdots,j_N}) - f_{\mathcal{Y}}(\theta_{j_1,\cdots,j_N})\right\vert. 
\end{eqnarray}

From Theorem~\ref{thm:3.2}, we have 
\begin{eqnarray}\label{eq:4.6-2}
\lefteqn{\left(\mathbb{E}_{\epsilon} \sup\limits_{\mathcal{X} \in  S_{\xi;(J_1,\cdots,J_N)}}\left\vert \sum\limits_{j_1=1,\cdots,j_N=1}^{J_1,\cdots,J_N}\epsilon_{i_1,\cdots,j_N}f^2_{\mathcal{X}}(\theta_{j_1,\cdots,j_N})\right\vert^p \right)^{1/p}} \nonumber \\
&\lesssim& \gamma_2\left(S_{\xi;(J_1,\cdots,J_N)},d_{\theta}\right) \sup\limits_{\mathcal{X} \in  S_{\xi;(J_1,\cdots,J_N)}}\left(\sum\limits_{j_1=1,\cdots,j_N=1}^{J_1,\cdots,J_N}f^2_{\mathcal{X}}(\theta_{j_1,\cdots,j_N})\right)^{1/2} \nonumber \\
& & +  \sup\limits_{\mathcal{X} \in  S_{\xi;(J_1,\cdots,J_N)}}\left(\mathbb{E}_{\epsilon}\left\vert \sum\limits_{j_1=1,\cdots,j_N=1}^{J_1,\cdots,J_N}\epsilon_{i_1,\cdots,j_N} f^2_{\mathcal{X}}(\theta_{j_1,\cdots,j_N})\right\vert^p \right)^{1/p},
\end{eqnarray}
where $\lesssim$ is a less than similar to comparison with a scalar constant difference. Since we have the following estimation from generic chaining:
\begin{eqnarray}
\gamma_2\left(S_{\xi;(J_1,\cdots,J_N)},d_{\theta}\right)&\lesssim& \int_0^{\infty}\left(\log N(S_{\xi;(J_1,\cdots,J_N)},d_{\theta},u)\right)^{1/2} du, 
\end{eqnarray}
where $N(S_{\xi;(J_1,\cdots,J_N)},d_{\theta},u)$ denote the covering number of the set $S_{\xi;(J_1,\cdots,J_N)}$, i.e., the smallest number of balls of radius $u$ in $(S_{\xi;(J_1,\cdots,J_N)}; d_{\theta})$ needed to cover $S_{\xi;(J_1,\cdots,J_N)}$. From~\cite{rudelson2008sparse}, we have the following bound about $\int_0^{\infty}\left(\log N(S_{\xi;(J_1,\cdots,J_N)},d_{\theta},u)\right)^{1/2} du$: 
\begin{eqnarray}\label{eq:4.6}
\int_0^{\infty}\left(\log N(S_{\xi;(J_1,\cdots,J_N)},d_{\theta},u)\right)^{1/2} du
&\lesssim& \Upsilon \sqrt{\frac{\xi}{\prod\limits_{i=1}^M I_i}} \log\left(\xi \sqrt{\log \prod\limits_{i=1}^M I_i}
\sqrt{\log \prod\limits_{j=1}^N J_j}\right). 
\end{eqnarray}

From Hoeffding's inequality, we have the follwing estimation for the term $\Scale[0.6]{\left(\mathbb{E}_{\epsilon}\left\vert \sum\limits_{j_1=1,\cdots,j_N=1}^{J_1,\cdots,J_N}\epsilon_{i_1,\cdots,j_N} f^2_{\mathcal{X}}(\theta_{j_1,\cdots,j_N})\right\vert^p \right)^{1/p}}$ in Eq.~\eqref{eq:4.6-2}:
\begin{eqnarray}\label{eq:4.6+1}
\lefteqn{\left(\mathbb{E}_{\epsilon}\left\vert \sum\limits_{j_1=1,\cdots,j_N=1}^{J_1,\cdots,J_N}\epsilon_{i_1,\cdots,j_N} f^2_{\mathcal{X}}(\theta_{j_1,\cdots,j_N})\right\vert^p \right)^{1/p}}\nonumber \\
&\leq&  \sqrt{p}\left(\sum\limits_{j_1=1,\cdots,j_N=1}^{J_1,\cdots,J_N} f^2_{\mathcal{X}}(\theta_{j_1,\cdots,j_N})\right)^{1/2}
\max\limits_{\substack{1\leq j_i \leq J_i,\\ \mbox{\tiny for $1\leq i \leq N$}}}\left\vert f_{\mathcal{X}}(\theta_{j_1,\cdots,j_N})\right\vert. 
\end{eqnarray}
From the $\xi$-sparsity of $\mathcal{X}$, we can have the following bound:
\begin{eqnarray}\label{eq:4.6+2}
\max\limits_{\substack{1\leq j_i \leq J_i,\\ \mbox{\tiny for $1\leq i \leq N$}}}\left\vert f_{\mathcal{X}}(\theta_{j_1,\cdots,j_N}) \right\vert \leq \Upsilon \sqrt{\frac{\xi}{\prod\limits_{i=1}^M I_i}}.
\end{eqnarray}
From estimations provided by Eqs.~\eqref{eq:4.6},~\eqref{eq:4.6+1} and~\eqref{eq:4.6+2}, we have 
\begin{eqnarray}
\lefteqn{\left(\mathbb{E}_{\epsilon}\left\vert \sum\limits_{j_1=1,\cdots,j_N=1}^{J_1,\cdots,J_N}\epsilon_{i_1,\cdots,j_N} f^2_{\mathcal{X}}(\theta_{j_1,\cdots,j_N})\right\vert^p \right)^{1/p}}\nonumber \\
&\leq&  \sup\limits_{\mathcal{X} \in  S_{\xi;(J_1,\cdots,J_N)}}\left(\sum\limits_{j_1=1,\cdots,j_N=1}^{J_1,\cdots,J_N}f^2_{\mathcal{X}}(\theta_{j_1,\cdots,j_N})\right)^{1/2}
\nonumber \\
&& \times \Upsilon \sqrt{\frac{\xi}{\prod\limits_{i=1}^M I_i}}\left(\sqrt{\log \prod\limits_{i=1}^M I_i}
\sqrt{\log \prod\limits_{j=1}^N J_j}\cdot \log\xi + \sqrt{p}\right).
\end{eqnarray}

From Eq.~\eqref{eq:4.5} and $L^p$ norm, we have
\begin{eqnarray}\label{eq:4.7}
\lefteqn{\left(\mathbb{E}\sup\limits_{\mathcal{X} \in  S_{\xi;(J_1,\cdots,J_N)}}\left\vert \sum\limits_{j_1=1,\cdots,j_N=1}^{J_1,\cdots,J_N}f^2_{\mathcal{X}}(\theta_{j_1,\cdots,j_N}) - \mathbb{E}f^2_{\mathcal{X}}(\theta_{j_1,\cdots,j_N}) \right\vert^p\right)^{1/p}} \nonumber \\
&\lesssim&  \left(\mathbb{E}\sup\limits_{\mathcal{X} \in  S_{\xi;(J_1,\cdots,J_N)}}\left\vert \sum\limits_{j_1=1,\cdots,j_N=1}^{J_1,\cdots,J_N}f^2_{\mathcal{X}}(\theta_{j_1,\cdots,j_N}) - \mathbb{E}f^2_{\mathcal{X}}(\theta_{j_1,\cdots,j_N}) \right\vert^p\right)^{1/2p} \nonumber \\
&& \times \Upsilon \sqrt{\frac{\xi}{\prod\limits_{i=1}^M I_i}}\left(\sqrt{\log \prod\limits_{i=1}^M I_i}
\sqrt{\log \prod\limits_{j=1}^N J_j}\cdot \log\xi + \sqrt{p}\right) \nonumber \\
&& + \Upsilon \sqrt{\frac{\xi}{\prod\limits_{i=1}^M I_i}}\left(\sqrt{\log \prod\limits_{i=1}^M I_i}
\sqrt{\log \prod\limits_{j=1}^N J_j}\cdot \log\xi  + \sqrt{p}\right).
\end{eqnarray}
Since Eq.~\eqref{eq:4.7} is the quadratic inequality related to $\left(\mathbb{E}\tau^p_{\xi}(\mathcal{U}_{\Pi})\right)^{1/2p}$, we can solve Eq.~\eqref{eq:4.7} with respect to $\left(\mathbb{E}\tau^p_{\xi}(\mathcal{U}_{\Pi})\right)^{1/2p}$ and obtain the following:
\begin{eqnarray}
\left(\mathbb{E}\tau^p_{\xi}(\mathcal{U}_{\Pi})\right)^{1/p} &\lesssim&  \Upsilon\sqrt{\frac{\xi}{\prod\limits_{i=1}^M I_i}}\left(\sqrt{\log \prod\limits_{i=1}^M I_i}
\sqrt{\log \prod\limits_{j=1}^N J_j}\cdot \log\xi \right) \nonumber \\
&& + \Upsilon^2\frac{\xi}{\prod\limits_{i=1}^M I_i} \left(\log \prod\limits_{i=1}^M I_i\right)\left(\log \prod\limits_{j=1}^N J_j \right) \cdot \log^2\xi \nonumber \\
&& +\Upsilon \sqrt{\frac{p\xi}{\prod\limits_{i=1}^M I_i}}+\Upsilon^2 \frac{p\xi}{\prod\limits_{i=1}^M I_i}.
\end{eqnarray}
From Lemma~\ref{lma:A1} and the condition provided by Eq.~\eqref{eq1:thm:4.1}, we have
\begin{eqnarray}
\mathrm{P}\left(\tau_{\xi} \geq \tau\right) &\leq&
\mathrm{P}\left(\tau_{\xi} \geq  G(\eta) \right) \leq \eta,
\end{eqnarray}
where $G(\eta)$ can be expressed as
\begin{eqnarray}
G(\eta)&=&\Upsilon\sqrt{\frac{\xi}{\prod\limits_{i=1}^M I_i}}\left(\sqrt{\log \prod\limits_{i=1}^M I_i}
\sqrt{\log \prod\limits_{j=1}^N J_j}\cdot \log\xi \right) \nonumber \\
&& + \Upsilon^2\frac{\xi}{\prod\limits_{i=1}^M I_i} \left(\log \prod\limits_{i=1}^M I_i\right)\left(\log \prod\limits_{j=1}^N J_j \right) \cdot \log^2\xi  +\Upsilon \sqrt{\frac{\log(\eta^{-1})\xi}{\prod\limits_{i=1}^M I_i}}+\Upsilon^2  \frac{\log(\eta^{-1})\xi}{\prod\limits_{i=1}^M I_i}.
\end{eqnarray}
Therefore, this theorem is proved.
$\hfill \Box$

\section{Application: Empirical Process}\label{sec:Application: Empirical Process} 

In this section, we will apply Theorem~\ref{thm:3.5} to study tail bounds for suprema of empirical processes. 
We require following lemma about Bernstein’s inequality for random tensors. 

\begin{lemma}\label{lma:Bernstein}
Given a finite sequence of independent Hermitian tensors $\{ \mathcal{X}_i  \in \mathbb{C}^{I_1 \times \cdots \times I_M  \times I_1 \times \cdots \times I_M } \}$ that satisfy
\begin{eqnarray}\label{eq1:thm:Subexponential Tensor Bernstein}
\mathbb{E} \mathcal{X}_i = 0 \mbox{~~and~~} \mathbb{E} (\mathcal{X}^p_i) \preceq \frac{p! \Upsilon^{p-2} }{2} \mathcal{A}_i^2,
\end{eqnarray}
where $p = 2,3,4,\cdots$, some positive contant $\Upsilon$. 

Define the total varaince $\sigma^2$ as: $\sigma^2 \define \left\Vert \frac{1}{n}\sum\limits_{i=1}^n \mathcal{A}_i^2 \right\Vert_{\mbox{\tiny spec}}$, where $\left\Vert\cdot \right\Vert_{\mbox{\tiny spec}}$ represents the spectral norm.
Then, we have following inequalities:
\begin{eqnarray}\label{eq2:thm:Subexponential Tensor Bernstein}
\mathrm{Pr} \left( \lambda_{\max}\left( \frac{1}{n}\sum\limits_{i=1}^{n} \mathcal{X}_i \right)\geq \sigma \sqrt{\frac{2}{n}} \sqrt{u} + \frac{\Upsilon}{n} u\right) \leq 2\left(\prod\limits_{i=1}^M I_i \right)\exp(-u).
\end{eqnarray}
\end{lemma}
\textbf{Proof:}
Slight modification from Theorem 13 in~\cite{chang2021convenient}.
$\hfill \Box$

We consider $n$ probability spaces $(\Omega_i, \mathrm{P}_i)$ for $i=1,2,\cdots,n$. Given a parameter set $T$ consisting of $n$-tuples $t=(t_1,\cdots, t_n)$, for every $t \in T$, we will have $n$-tuple $\mathcal{X}_t = \left(\mathcal{X}_{t_1}, \cdots,\mathcal{X}_{t_n}\right)$ of random Hermitian tensors satisfying those conditions required by Lemma~\ref{lma:Bernstein}. We consider the following exmpirical tensor process 
\begin{eqnarray}\label{eq:Et def}
\mathcal{E}_t &=& \frac{1}{n}\sum\limits_{i=1}^n\left(\mathcal{X}_{t_i} - \mathbb{E}\mathcal{X}_{t_i}\right). 
\end{eqnarray}
From Lemma~\ref{lma:Bernstein}, we have the tensor process $(\mathcal{E}_t)_{t \in T}$ has a mixed tail with respect to the following metrics 
\begin{eqnarray}
d_1(s,t)&=&\max\limits_{1 \leq i \leq n} \lambda_{\max}(\mathcal{X}_{t_i} - \mathcal{X}_{s_i}),\nonumber \\
d_2(s,t)&=&\left(\frac{1}{n}\sum\limits_{i=1}^n \lambda_{\max}(\mathcal{X}_{t_i} - \mathcal{X}_{s_i})^2\right)^{1/2}.
\end{eqnarray}

\begin{theorem}\label{thm:cor4.6}
Let $\mathcal{E}_t$ be defined as Eq.~\eqref{eq:Et def} and satisfy the following:
\begin{eqnarray}\label{eq1:thm:cor4.6}
\mathcal{X}_{t_i} - \mathbb{E}\mathcal{X}_{t_i}\preceq \frac{p! \Upsilon^{p-2} }{2} \mathcal{A}_i^2,
\end{eqnarray}
then, for any $1 \leq p < \infty$, we have
\begin{eqnarray}
\left(\mathbb{E}\sup\limits_{t \in T} \left\Vert\mathcal{E}_t\right\Vert^p_{\mbox{\tiny spec}} \right)^{1/p}
&\lesssim&\left(\frac{1}{\sqrt{n}}\gamma_2(T,d_2)+\frac{1}{n}\gamma_1(T,d_1)\right)+\sqrt{p}\frac{\sigma}{\sqrt{n}} + p\frac{\Upsilon}{n},
\end{eqnarray}
where $\sigma^2 \define \left\Vert \frac{1}{n}\sum\limits_{i=1}^n \mathcal{A}_i^2 \right\Vert_{\mbox{\tiny spec}}$. 

Moreover, there exist contants $C_2, C_3 >0$ such that for any $u \geq 1$, we have
\begin{eqnarray}
\mathrm{P}\left(\sup\limits_{t \in T} \left\Vert\mathcal{E}_t\right\Vert_{\mbox{\tiny spec}}\geq
C_2\left(\frac{1}{\sqrt{n}}\gamma_2(T,d_2)+\frac{1}{n}\gamma_1(T,d_1)\right) +C_3\left(\frac{\sigma}{\sqrt{n}}\sqrt{u} + \frac{\Upsilon}{n}u\right)\right) \leq e^{-u}.
\end{eqnarray}
\end{theorem}
\textbf{Proof:}
Apply Theorem~\ref{thm:3.5} with Lemma~\ref{lma:Bernstein} by setting $\left\Vert \cdot \right\Vert_{\alpha}$ as 
$\left\Vert\cdot\right\Vert_{\mbox{\tiny spec}}$.
$\hfill \Box$

\bibliographystyle{IEEETran}
\bibliography{TensorSum_GC_Bib}

\end{document}